\documentclass[12pt]{amsart}
\usepackage{amsfonts, amsmath, amsthm}
\usepackage{graphicx}
\usepackage{times}\usepackage{geometry}

\newtheorem{pro}{Proposition}
\newtheorem{thm}[pro]{Theorem}
\newtheorem{lem}[pro]{Lemma}

\newtheorem{cnj}[pro]{Conjecture}

\newtheorem{cor}[pro]{Corollary}

\theoremstyle{definition}

\theoremstyle{remark}
\newtheorem*{rmk}{Remark}

\newcommand{\s}{\Sigma}

\title[Knots with $g(E(K)) = 2$ and $g(E(K \# K \# K)) = 6$]{Knots with
  $g(E(K)) = 2$ and $g(E(K \# K \# K)) = 6$ \\[3pt] and Morimoto's Conjecture} 

\date{\today} \address{Department of Mathematics, Nara Women's University
Kitauoya Nishimachi, Nara 630-8506, Japan} \address{Department of mathematical
Sciences, University of Arkansas, Fayetteville, AR 72701}
\email{tsuyoshi@cc.nara-wu.ac.jp} \email{yoav@uark.edu} \author{Tsuyoshi
Kobayashi} \author{Yo'av Rieck}

\begin{document}


\date{\today}%
\dedicatory{Dedicated to the memory of Yves Mathieu \\
and to Michel Domergue on the occasion of his 60th birthday}
\begin{abstract}
We show that there exist knots $K \subset S^3$ with $g(E(K)) = 2$ and $g(E(K
\# K \# K)) = 6$.  Together with \cite[Theorem~1.5]{crelle}, this proves
existence of counterexamples to Morimoto's Conjecture \cite{mc}. 
This is a special case of \cite{big}.
\end{abstract}
\maketitle

Let $K_i$ ($i=1,2$) be knots in the 3-sphere $S^3$, and let $K_1 \# K_2$ be
their connected sum.  We use the notation $t(\cdot)$, $E(\cdot)$, and
$g(\cdot)$ to denote tunnel number, exterior, and Heegaard genus respectively
(we follow the definitions and notations given in \cite{cag}).  It is well
known that the union of a tunnel system for $K_1$, a tunnel system for $K_2$,
and a tunnel on a decomposing annulus for $K_1 \# K_2$ forms a tunnel system
for $K_1 \# K_2$.  Therefore:
$$t(K_1 \# K_2) \leq t(K_1) + t(K_2) +1.$$
Since (for any knot $K$) $t(K) = g(E(K)) - 1$ this gives:
\begin{equation}
\label{eq:upper-bound}  
g(E(K_1 \# K_2)) \leq g(E(K_1)) + g(E(K_2)).
\end{equation}
We say that a knot $K$ in a closed orientable manifold $M$ admits a $(g,n)$ position
if there exists a genus $g$ 
Heegaard surface $\s \subset M$, separating $M$ into the handlebodies
$H_1$ and $H_2$, so that $H_i \cap K$ ($i=1,2$) consists of $n$ arcs that are
simultaneously parallel into $\partial H_i$.  It is known
\cite[Proposition~1.3]{mc} that if $K_i$ ($i=1$ or 2) admits a $(t(K_i),1)$ position
then equality does not hold:
$$g(E(K_1 \# K_2)) < g(E(K_1)) + g(E(K_2)).$$
Morimoto proved that if $K_1$ and $K_2$ are m-small knots then the converse
holds, and conjectured that this holds in general \cite[Conjecture~1.5]{mc}:
\begin{cnj}[Morimoto's Conjecture]
\label{cnj:mc}
Given knots $K_1,\ K_2 \subset S^3$, $g(E(K_1 \# K_2)) < g(E(K_1)) +
g(E(K_2))$ if and only if for $i=1$ or $i=2$, $K_i$ admits a $(t(K_i),1)$
position.  
\end{cnj}

We denote the connected sum of $n$ copies of $K$ by $nK$.  We prove:

\begin{thm}
\label{thm:additive}
There exists infinitely many knots $K \subset S^3$ with $g(E(K)) = 2$ and
$g(E(3K)) = 6$.  
\end{thm}

\begin{rmk}
This is a special case of \cite[Theorem~1.4]{big}.  By specializing we obtain an
easy and accessible argument that can be used as an introduction to the main
ideas of \cite{big}.
\end{rmk}

As in \cite{big} Theorem~\ref{thm:additive} implies:

\begin{cor}
\label{cor:mc}
There exists a counterexample to Morimoto's Conjecture, specifically, there
exist knots $K_1, \ K_2 \subset S^3$ so that $K_i$ does not admit a
$(t(K_i),1)$ position ($i=1,2$), and (for some integer $m$) $g(E(K_1)) = 4$,
$g(E(K_2)) = 2(m-2)$, and $g(E(K_1 \# K_2)) < 2m$.
\end{cor}

\begin{proof}[Proof of Corollary~\ref{cor:mc}]
See the proof of \cite[Corollary~1.8]{big}.
\end{proof}

We note that $K_1$ and $K_2$ are composite knots.  This leads Moriah
\cite[Conjecture~7.14]{moriah-s} to conjecture that if $K_1$ and $K_2$ are
prime then Conjecture~\ref{cnj:mc} holds.

\section{The proof.}

Let $X$ be the exterior of a knot $K$ in a closed orientable manifold.  For an
integer $c \geq 0$ let $X^{(c)}$ denote the manifold obtained by drilling $c$
curves out of $X$ that are simultaneously parallel to meridians of $K$.  The
following is \cite[Proposition~2.2]{big}, where the proof can be found.  Note
the relation to \cite[Theorem~3.8]{schsch}.

\begin{pro}
\label{pro:bridge}
Let $X$, $X^{(c)}$ be as above and $g \geq 0$ an integer.  Suppose $X^{(c)}$
admits a strongly irreducible Heegaard surface of genus $g$. Then one of the
following holds:
\begin{enumerate}
\item $X$ admits an essential surface $S$ with $\chi(S) \geq 4 - 2g$.
\item For some $b$, $c \leq b \leq g$, $K$ admits $(g-b,b)$ position.
\end{enumerate}
\end{pro}

Given an integer $d > 0$, Johnson and Thompson \cite{jt} and Minsky, Moriah
and Schleimer \cite{mms} construct infinitely many knots $K \subset S^3$ so
that $E(K)$ admits a genus 2 Heegaard splitting of distance more than $d$ (in
the sense of the curve complex \cite{hempel}).  (Note that \cite{mms}
is more general.)  Fix such a knot $K$ for $d=10$.
The two properties of $K$ we will need
are described in the lemmas below: 

\begin{lem}
\label{lem:essential}
$X$ does not admit an essential surface $S$ with $\chi(S) \geq -8$.
\end{lem}

\begin{proof}
This follows directly from \cite[Theorem~31.]{scharlemann-2004}.  
\end{proof}

\begin{lem}
\label{lem:bridge}  
$K$ does not admit a $(0,3)$ or a $(1,2)$ position.
\end{lem}

\begin{proof}
Assume, for a contradiction, $K$ admits a $(0,3)$ or a $(1,2)$ position.
By \cite[Theorem~1]{jt}, if $K$ admits a $(p,q)$ position (for some $p$,
$q$) then either $K$ is isotopic into a genus $p$ Heegaard surface, or
the distance of any Heegaard splitting of $X$ is at most $2(p+q)$.
Since $K$ is not a trivial knot or a torus knot, the former cannot happen.
(Note that, by \cite{schultens2} we see that the distance of each Heegaard 
splitting of the exterior of any torus knot is at most 2.) 
On the other hand, if the latter holds, then the distance of any
Heegaard splitting of $X$ should be at most $6$ contradicting our choice of $K$.

\end{proof}

For integers $n \geq 1$ and $c \geq 0$ we denote the exterior of
$nK$ by $X(n)$, and the manifold obtained 
by drilling $c$ curves out of $X(n)$ that are simultaneously parallel to
meridians of $nK$ by $X(n)^{(c)}$.  

Thus we obtain $X(n)^{(c)}$ by drilling a curve $\gamma_n \subset
X(n)^{(c-1)}$ that is parallel to $\partial X$, and in particular, $\gamma_n$
can be isotoped onto any Heegaard surface
of $X^{(c-1)}$.  This is described in \cite{rieck} by saying that $X^{(c-1)}$
is obtained from $X^{(c)}$ by a {\it good} Dehn filling.  For good Dehn
fillings \cite{rieck} shows (see the proof of Theorem~5.1 of \cite{cag} for
details): 

\begin{lem}
\label{lem:dehn}
Either $g(X(n)^{(c)}) = g(X(n)^{(c-1)})$ or $g(X(n)^{(c)} = g(X(n)^{(c-1)})+1.$
\end{lem}

\begin{lem}
\label{lem:g(x)1}  
$g(X^{(1)}) = 3$.
\end{lem}

\begin{proof}
Since $g(X) = 2$, by Lemma~\ref{lem:dehn} $g(X^{(1)}) = 2$ or $g(X^{(1)}) = 3$.
Assume for a contradiction that $g(X^{(1)}) = 2$ and let $\s^{(1)} \subset
X^{(1)}$ be a minimal genus Heegaard surface. 

\medskip

\noindent{Claim.} $\s^{(1)}$ is strongly irreducible.

\begin{proof}
Suppose $\s^{(1)}$ weakly reduces.  Then by Casson and Gordon
\cite{casson-gordon} (see \cite[Theorem~1.1]{sedgwick} 
for a relative
version) an appropriately chosen weak reduction yields an 
essential surface $S$ with $\chi(S) \geq \chi(\s^{(1)}) + 4 = 2$.  Since
$X^{(1)}$ does not admit an essential sphere this is a contradiction, proving
the claim.
\end{proof}

Thus we may assume $\s^{(1)}$ is strongly irreducible.  By
Proposition~\ref{pro:bridge}, either $X$ admits an essential surface $S$ with
$\chi(S) \geq 4 - 2g(\s^{(1)}) = 0$ or $K$ admits a $(2-b,b)$ position,
with $1 \leq b \leq 2$.  The former contradicts Lemma~\ref{lem:essential}.  For
the latter, we get a $(1,1)$ position (for $b=1$) or a $(0,2)$ position (for
$b=2$).  Both contradict Lemma~\ref{lem:bridge}.  
\end{proof}

\begin{lem}
\label{lem:g(x)2}  
$g(X^{(2)}) = 4$.
\end{lem}

\begin{proof}
Since $g(X^{(1)}) = 3$, by Lemma~\ref{lem:dehn} $g(X^{(2)}) = 3$ or
$g(X^{(2)}) = 4$.  Assume for a contradiction that $g(X^{(2)}) = 3$  and let
$\s^{(2)} \subset X^{(2)}$ be a minimal genus Heegaard surface.  

\medskip

\noindent{Claim.} $\s^{(2)}$ is strongly irreducible.

\begin{proof}
Suppose $\s^{(2)}$ weakly reduces.  Then by Casson and Gordon
\cite{casson-gordon} (see \cite{sedgwick} for a relative version) an
appropriately chosen weak reduction yields an  essential surface $S$ with
$\chi(S) \geq \chi(\s^{(2)}) + 4 = 0$.  Since $X^{(2)}$ does not admit an
essential sphere, this surface must be a collection of tori; let $F$ be one of
these tori.  By \cite[Proposition~2.13]{cag}, $\s^{(2)}$ weakly reduces to $F$.

Note that $X^{(2)}$ admits an essential torus $T$ giving the decomposition
$X^{(2)} = X' \cup_T Q^{(2)}$, where $Q^{(2)}$ is homeomorphic to an annulus
with two holes cross ~$S^1$ and $X' \cong X$.

Since $F$ and $T$ are incompressible, we may suppose that each component of $F
\cap T$ is a simple closed curve which is essential in both $F$ and $T$.
Minimize $|F \cap T|$ under this constraint.  We claim that $F \cap T =
\emptyset$.  Assume for a contradiction $F \cap T \neq \emptyset$.  Then any
component of $F \cap X'$ is an essential annulus; by
Lemma~\ref{lem:essential},  $X'$ does not admit essential annuli.

Thus we may assume $F \subset X'$ or $F \subset Q^{(2)}$.  If $F \subset X'$
and not parallel to $T$ then
$X \cong X'$ is toroidal, contradicting Lemma~\ref{lem:essential}.  If $F$ is
parallel to $T$ we isotope it into $Q^{(2)}$. 

Thus we may assume $F \subset Q^{(2)}$. By \cite[VI.34]{jaco} $F$ 
is a vertical torus in $Q^{(2)}$.  Assume first that
$F$ is isotopic to a component of $\partial Q^{(2)}$.  Since $F$ was obtained
by weakly reducing a minimal genus Heegaard surface for $X^{(2)}$, by
\cite[Theorem~1.1]{sedgwick} $F$ is not peripheral, {\it i.e.}, 
$F$ is not isotopic to a component of $\partial X^{(2)}$.  Hence $F$ is
isotopic to  $T$ and $X^{(2)} = X' \cup_F Q^{(2)}$.  Note that by
\cite{schultens} $g(Q^{(2)}) = 3$, and since $X \cong X'$, $g(X') = 2$.  Since
$F$ was obtained by weakly reducing a {\it minimal genus} Heegaard surface, 
\cite[Proposition~2.9]{cag} (see also \cite[Remark~2.7]{schultens}) gives:  
$$g(X^{(2)}) = g(Q^{(2)}) + g(X') - g(F) = 3 + 2 - 1 = 4.$$ This contradicts
our assumption that $g(X^{(2)}) = 3$.

Next assume that $F$ is not isotopic to a component of $\partial Q^{(2)}$.
Then $F$ is isotopic to a vertical torus giving the decomposition $X^{(2)} =
X_1 \cup_F D(2)$, where $X_1$ is homeomorphic to $X^{(1)}$ and $D(2)$ is
homeomorphic to a twice punctured disk cross $S^1$.  By Lemma~\ref{lem:g(x)1}
$g(X_1) = 3$ and by \cite{schultens} $g(D(2)) = 2$.  We get:
$$g(X^{(2)}) = g(X_1) + g(D(2)) - g(F) = 3 + 2 - 1 = 4.$$ This contradicts
our assumption that $g(X^{(2)}) = 3$.

This contradiction proves the claim.
\end{proof}

Thus we may assume $\s^{(2)}$ is strongly irreducible.  By
Proposition~\ref{pro:bridge}, either $X$ admits an essential surface $S$ with
$\chi(S) \geq 4 - 2g(\s^{(2)}) = -2$ or $K$ admits a $(g(\s^{(2)})-b,b) =
(3-b,b)$ position, with $2 \leq b \leq 3$.  The former contradict
Lemma~\ref{lem:essential}.  For the 
latter, we get a $(1,2)$ position (for $b=2$) or a $(0,3)$ position (for
$b=3$).  Both contradict Lemma~\ref{lem:bridge}.
\end{proof}

\begin{lem}
\label{lem:g(x(2))}  
$g(X(2)) = 4$.
\end{lem}

\begin{proof}
By Inequality~(\ref{eq:upper-bound}) $g(X(2)) \leq 4$.  Therefore by the Swallow
Follow Torus Theorem \cite[Theorem~4.1]{cag} and Lemma~\ref{lem:essential} any
minimal genus Heegaard surface for $X(2)$ weakly reduces to a swallow follow
torus $F$, giving the decomposition:
$X(2) = X^{(1)} \cup_F X.$ By \cite[Proposition~2.9]{cag} and
Lemma~\ref{lem:g(x)1}, $g(X(2)) = g(X^{(1)}) 
  + g(X) - g(F) = 3 + 2 - 1 = 4$.  
\end{proof}

\begin{lem}
\label{lem:g(X(2)1)}  
$g(X(2)^{(1)}) = 5$.
\end{lem}

\begin{proof}
By Lemmas~\ref{lem:dehn} and \ref{lem:g(x(2))}, $g(X(2)^{(1)}) = 4$ or
$g(X(2)^{(1)}) = 5$.  Assume for a contradiction that $g(X(2)^{(1)}) = 4$.
By the Swallow Torus Theorem \cite[Theorem~4.2]{cag} and
Lemma~\ref{lem:essential} any minimal genus Heegaard surface for $X(2)^{(1)}$
weakly reduces to a swallow follow torus $F$ giving one of the following
decompositions: 
\begin{enumerate}
\item $X(2)^{(1)} = X(2) \cup_F Q^{(1)}$, where $Q^{(1)}$ is
homeomorphic to an annulus with one hole cross $S^1$.
\item $X(2)^{(1)} = X^{(1)} \cup_F X^{(1)}$.
\item $X(2)^{(1)} = X^{(2)} \cup_F X$.
\end{enumerate}
By \cite{schultens} $g(Q^{(1)}) = 2$; the genera of all other manifolds are given
in the lemmas above.  By amalgamation \cite[Proposition~2.9]{cag} we get:
\begin{enumerate}
\item $g(X(2)^{(1)}) = g(X(2)) + g(Q^{(1)}) - g(F) = 4 + 2 -1 =5$.
\item $g(X(2)^{(1)}) = g(X^{(1)}) + g(X^{(1)}) - g(F) = 3 + 3  - 1 = 5$.
\item $g(X(2)^{(1)}) = g(X^{(2)}) + g(X)- g(F) = 4 + 2 - 1 = 5$.
\end{enumerate}
\end{proof}

\begin{proof}[Proof of Theorem~\ref{thm:additive}]
By Inequality~(\ref{eq:upper-bound}), $g(X(3)) \leq 6$.  Therefore, by the 
Swallow Follow Torus Theorem \cite[Theorem~4.2]{cag} and
Lemma~\ref{lem:essential} any minimal genus Heegaard surface for $X(3)$ weakly
reduces to a swallow follow torus $F$ giving one of the following
decompositions:
\begin{enumerate}
\item $X(3) = X^{(1)} \cup_F X(2)$.
\item $X(3) = X(2)^{(1)} \cup_F X$.
\end{enumerate}
The genera of the manifolds are given in the lemmas above.  By amalgamation
\cite[Proposition~2.9]{cag} we get: 
\begin{enumerate}
\item $g(X(3)) = g(X^{(1)}) + g(X(2)) - g(F) = 3 + 4 -1 = 6$.
\item $g(X(3)) = g(X(2)^{(1)}) + g(X) - g(F) = 5 + 2 - 1 = 6$.
\end{enumerate}
  
This completes the proof of Theorem~\ref{thm:additive}.
\end{proof}



\begin{thebibliography}{10}

\bibitem{casson-gordon}
A.~J. Casson and C.~McA. Gordon.
\newblock Reducing {H}eegaard splittings.
\newblock {\em Topology Appl.}, 27(3):275--283, 1987.

\bibitem{hempel}
John Hempel.
\newblock 3-manifolds as viewed from the curve complex.
\newblock {\em Topology}, 40(3):631--657, 2001.

\bibitem{jaco}
William Jaco.
\newblock {\em Lectures on three-manifold topology}, volume~43 of {\em CBMS
  Regional Conference Series in Mathematics}.
\newblock American Mathematical Society, Providence, R.I., 1980.

\bibitem{jt}
Jesse Johnson and Abigail Thompson.
\newblock On tunnel number one knots which are not (1,n), 2006.
\newblock available at http://arxiv.org/abs/math.GT/0606226.

\bibitem{crelle}
Tsuyoshi Kobayashi and Yo'av Rieck.
\newblock On the growth rate of the tunnel number of knots.
\newblock {\em J. Reine Angew. Math.}, 592:63--78, 2006.

\bibitem{big}
Tsuyoshi Kobayashi and Yo'av Rieck.
\newblock Knot exteriors with additive {H}eegaard genus and {M}orimto's
{C}onjecture. 
\newblock available at http://arxiv.org/abs/math.GT/0701765, 2007.

\bibitem{cag}
Tsuyoshi Kobayashi and Yo'av Rieck.
\newblock Heegaard genus of the connected sum of m-small knots.
\newblock {\em Communications in Analysis and Geometry}, to appear 2006.

\bibitem{mms}
Yair Minsky, Yoav Moriah, and Saul Schleimer.
\newblock High distance knots, 2006.
\newblock available at http://arxiv.org/abs/math.GT/0607265.

\bibitem{moriah-s}
Yoav Moriah.
\newblock Heegaard splittings of knot exteriors.
\newblock http://arxiv.org/abs/math.GT/0608137, 2006.

\bibitem{mc}
Kanji Morimoto.
\newblock On the super additivity of tunnel number of knots.
\newblock {\em Math. Ann.}, 317(3):489--508, 2000.

\bibitem{rieck}
Yo'av Rieck.
\newblock Heegaard structures of manifolds in the {D}ehn filling space.
\newblock {\em Topology}, 39(3):619--641, 2000.

\bibitem{scharlemann-2004}
Martin Scharlemann.
\newblock Proximity in the curve complex: boundary reduction and bicompressible
  surfaces, 2004.
\newblock http://arxiv.org/abs/math.GT/0410278.

\bibitem{schsch}
Martin Scharlemann and Jennifer Schultens.
\newblock Comparing {H}eegaard and {JSJ} structures of orientable 3-manifolds.
\newblock {\em Trans. Amer. Math. Soc.}, 353(2):557--584 (electronic), 2001.

\bibitem{schultens}
Jennifer Schultens.
\newblock The classification of {H}eegaard splittings for (compact orientable
  surface){$\,\times\, S\sp 1$}.
\newblock {\em Proc. London Math. Soc. (3)}, 67(2):425--448, 1993.


\bibitem{schultens2}
Jennifer Schultens.
\newblock Heegaard splittings of Seifert fibered spaces with boundary.
\newblock {\em Trans. Amer. Math. Soc.}, 347:2533--2552, 1995.

\bibitem{sedgwick}
Eric Sedgwick.
\newblock Genus two 3-manifolds are built from handle number one pieces.
\newblock {\em Algebr. Geom. Topol.}, 1:763--790 (electronic), 2001.

\end{thebibliography}

\end{document}